# EXACT INTEGRATION SCHEME FOR SIX-NODE WEDGE ELEMENT MASS MATRIX.

E Hanukah

Faculty of Mechanical Engineering, Technion – Israel Institute of Technology, Haifa 32000, Israel
Corresponding author Email: eliezerh@tx.technion.ac.il

**Abstract**

Currently, mass matrices are computed by means of sufficiently accurate numerical integration schemes. Two-point and nine-point (Gauss) quadrature remain frequently used. We derive an exact, easy to implement integration rule for six-node wedge element mass matrices based on seven points only. Both consistent and lumped mass matrices have been considered. New metric (jacobian determinant) interpolation accompanied by analytical integration permits computing effort reduction next to accuracy increase of integration rule. In addition, one and four point mass matrices integration schemes have been proposed. Accuracy superiority over equivalent schemes is established.

**Key words**: numerical integration, quadrature, prism finite element, pentahedron, closed-form, symbolic computational mechanics, semi-analytical, mass matrix.

## 1. Introduction

Mass coefficients, internal forces, stiffness matrix, all require integration in the element domain, which is most commonly obtained with the help of numerical integration schemes [1]. Several studies exploit the idea of closed-form integration for stiffness matrixes [2-9], significant time savings has been established. Furthermore, hierarchical semi-analytical displacement based approach is used to model three dimensional finite bodies e.g. [10-12] yielding new analytical solutions.

In present study we follow the basic guidelines presented in [13, 14]. The metric is approximated using interpolations of different order. Zero order approximation requires metric evaluation at one point, usually the centroid; first order interpolation involves four points and linear interpolation function, etc. Analytical integration is performed, whereas convenient coefficient matrices definition allows familiar representation of the resulting integration rule.

We consider constant metric (CM) approximation; linear metric (LM) approximation; special interpolation, based on seven points, which exactly (EX) represent the metric. Analytical integration and coefficient matrices definition permits familiar representation, namely, CM mass matrix is equivalent to one point numerical scheme, LM corresponds to four point rule, and EX is similar to seven point scheme.

New schemes are exact for fine mesh (constant metric element, triangular bases are parallel to each other), hence, preliminary numerical study examines coarse mesh accuracy of the schemes and compares to currently used two and nine points rules. Accuracy superiority for entire coarseness range has been established.

The outline of the paper as follows. Section 2 recalls important details of six node solid wedge element concluding with mass matrix formulation. Natural coordinates, shape functions, metric and jacobian matrix, numerical and analytical integration in the element domain are recalled. Section 3 starts with jacobian matrix and the metric representation. Special CM, LM, and EX metric interpolations are proposed. Section 4 presents mass matrix details following analytical integration. Coefficient matrices are explicitly presented. Section 5 contains preliminary numerical accuracy study including comparison to equivalent numerical integration schemes. Coarse mesh has been considered. Section 6 summarizes and records our conclusions.

## 2. Background

Initial nodal locations of standard six node wedge element (e.g. [15] pp.73) denoted by $\mathbf{X}_i$ $(i=1,..,6)$. $X_{mi}$ $(m=1,2,3, i=1,..,6)$ stand for nodal components in terms of global Cartesian coordinates system $\mathbf{X}_i = X_{mi}\mathbf{e}_m$ $(m=1,2,3, i=1,..,6)$, summation convention on repeated index is implied. Here and throughout the text, bold symbols traditionally denote vector or tensor quantities. Local convected coordinate system (i.e. natural coordinates) $\{\xi,\eta,\zeta\}$ see fig. 1 admits

$$0 \leq \xi \leq 1-\eta \ , \ 0 \leq \eta \ , \ \zeta \leq 1 \qquad (1)$$

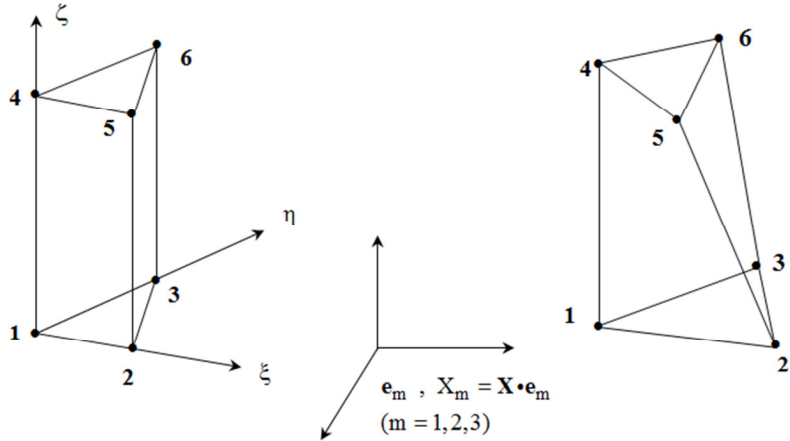

Figure 1: Showing the Global coordinate system $\mathbf{e}_m$ $(m=1,2,3)$, natural coordinates $\{\xi,\eta,\zeta\}$, nodes numbering, straight-sided element on the left and general (coarse) mesh on the right.

Shape functions $\varphi^i(\xi,\eta,\zeta)(i=1,..,6)$ in terms of natural coordinates are given by

$$\varphi^1 = \frac{1}{2}(1-\xi-\eta)(1-\zeta) \;,\; \varphi^2 = \frac{\xi}{2}(1-\zeta) \;,\; \varphi^3 = \frac{\eta}{2}(1-\zeta)$$
$$\varphi^4 = \frac{1}{2}(1-\xi-\eta)(1+\zeta) \;,\; \varphi^5 = \frac{\xi}{2}(1+\zeta) \;,\; \varphi^6 = \frac{\eta}{2}(1+\zeta) \tag{2}$$

Location of material point X inside the element domain (1) is denoted by **X** and given by

$$\mathbf{X} = \varphi^i \mathbf{X}_i \;\; (i=1,..,6) \tag{3}$$

Homogeneous initial configuration $\rho_0 = \text{const}$ is considered. Extension for linearly varying initial density $\rho_0 = (1/2-\xi-\eta-\zeta/2)\rho_1 + \xi\rho_2 + \eta\rho_3 + (\zeta/2+1/2)\rho_4$ or isoparametric initial density $\rho_0 = \varphi^i \rho_i$ (i=1,..,6) where $\rho_i$ denote initial nodal densities, follows the similar stages. Metric or jacobian determinant of global-local coordinate's transformation J is given by

$$J = \mathbf{X}_{,1} \times \mathbf{X}_{,2} \cdot \mathbf{X}_{,3} = \begin{vmatrix} (\mathbf{X}\cdot\mathbf{e}_1)_{,1} & (\mathbf{X}\cdot\mathbf{e}_1)_{,2} & (\mathbf{X}\cdot\mathbf{e}_1)_{,3} \\ (\mathbf{X}\cdot\mathbf{e}_2)_{,1} & (\mathbf{X}\cdot\mathbf{e}_2)_{,2} & (\mathbf{X}\cdot\mathbf{e}_2)_{,3} \\ (\mathbf{X}\cdot\mathbf{e}_3)_{,1} & (\mathbf{X}\cdot\mathbf{e}_3)_{,2} & (\mathbf{X}\cdot\mathbf{e}_3)_{,3} \end{vmatrix} > 0 \tag{4}$$

$$J_{mn}(\xi,\eta,\zeta;X_{ki}) = (\mathbf{X}\cdot\mathbf{e}_m)_{,n} \;,\; (i=1,..,6, m,n,k=1,2,3)$$

Where (×) and (·) stand for vector cross and scalar products, |·| stand for determinant operator, comma denotes partial differentiation with respect to natural coordinates. Here and throughout the text, determinant of jacobian 3x3 matrixes is computed by

$$J = J_{11}J_{22}J_{33} - J_{11}J_{23}J_{32} - J_{31}J_{22}J_{13} - J_{21}J_{12}J_{33} + J_{21}J_{32}J_{13} + J_{31}J_{12}J_{23} \tag{5}$$

The above is consistent with standard determinant definition. Differential volume element dV and initial volume V are defined by

$$dV = J(\xi,\eta,\zeta;X_{mi})d\xi d\eta d\zeta \;\; (m=1,2,3, i=1,..,6) \;,\; V = \int_V dV \tag{6}$$

Isoparametric formulation (e.g.[16] pp.104) for mass conserving element yield consistent, symmetric, positive definite mass matrix

$$M^{ij} = \int_V \rho_0 \phi^i \phi^j dV = \rho_0 \int_0^{+1}\int_0^{+1}\int_0^{1-\eta} \phi^i \phi^j J d\xi d\eta d\zeta \;,\; M^{ij} = M^{ij} \;,\; (i,j=1,..,6) \tag{7}$$

Widely used (e.g. [17-19]) diagonal or lumped mass matrix is given by

$$M_{diag}^{ii} = \int_V \rho_0 \phi^i dV = \rho_0 \int_0^{+1}\int_0^{+1}\int_0^{1-\eta} \phi^i J d\xi d\eta d\zeta \;,\; M_{diag}^{ij} = \begin{cases} M_{diag}^{ii} & i=j \\ 0 & i \neq j \end{cases} \;,\; (i,j=1,..,6) \tag{8}$$

Numerical integration of (7) result in (e.g. [16] pp.126)

$$\int_V \rho_0 \phi^i \phi^j J d\xi d\eta d\zeta \approx \rho_0 \sum_{p=1}^{n_p} w_p \phi^i(\xi_p,\eta_p,\zeta_p) \phi^j(\xi_p,\eta_p,\zeta_p) J(\xi_p,\eta_p,\zeta_p;X_{mi}) \quad (9)$$

$$(m = 1,2,3, i = 1,..,6)$$

Where $n_p$ stand for number of integration (Gauss) points, $w_p$ denotes weights and $\xi_p, \eta_p, \zeta_p$ are coordinates of integration points. Two and nine points schemes are broadly accepted [15] pp.80. For later convenience, the next definitions are considered

$$\hat{J}_p(X_{mi}) = J(\xi_p,\eta_p,\zeta_p;X_{mi}) \;,\; \hat{M}_p^{ij} = w_p \phi^i(\xi_p,\eta_p,\zeta_p) \phi^j(\xi_p,\eta_p,\zeta_p) \quad (10)$$

$$(m=1,2,3, i,j=1,..,6, p=1,..,n_p)$$

Such that

$$n_p = 2 \;,\; M^{ij} = \rho_0(\hat{J}_1 \hat{M}_1^{ij} + \hat{J}_2 \hat{M}_2^{ij}) \;,\; (i,j=1,..,6)$$
$$n_p = 9 \;,\; M^{ij} = \rho_0(\hat{J}_1 \hat{M}_1^{ij} + ... + \hat{J}_4 \hat{M}_4^{ij} + ... + \rho_0 \hat{J}_9 \hat{M}_9^{ij}) \quad (11)$$

Herein we'll show that by deriving successful interpolation for the metric, and performing analytical integration, one develops new coefficient matrices $\hat{M}_p^{ij}$ $(i,j=1,..,6, p=1,..,n_p)$, and selects new "integration points" resulting in beneficial schemes, i.e. fewer integration points yet superior accuracy. Symbolic manipulations have been performed using computer algebra system MAPLE$^{TM}$.

### 3. Metric interpolation.

Jacobian matrix (4) is represented by

$$J_{mn} = J_{mn}^0 + \xi J_{mn}^1 + \eta J_{mn}^2 + \zeta J_{mn}^3 \;,\; (m,n=1,2,3) \quad (12)$$

Matrices $J_{mn}^k$ $(m,n=1,2,3, k=0,..,3)$ are functions of nodal components $X_{ki}$ $(k=1,2,3, i=1,..,6)$

$$J_{mn}^0 = \frac{1}{2} \begin{pmatrix} -X_{14}-X_{11}+X_{15}+X_{12} & X_{16}-X_{14}+X_{13}-X_{11} & -X_{11}+X_{14} \\ -X_{24}-X_{21}+X_{25}+X_{22} & X_{26}-X_{24}+X_{23}-X_{21} & -X_{21}+X_{24} \\ -X_{34}-X_{31}+X_{35}+X_{32} & X_{36}-X_{34}+X_{33}-X_{31} & -X_{31}+X_{34} \end{pmatrix} \quad (13)$$

$$J_{mn}^1 = \frac{1}{2} \begin{pmatrix} 0 & 0 & X_{11}-X_{12}-X_{14}+X_{15} \\ 0 & 0 & X_{21}-X_{22}-X_{24}+X_{25} \\ 0 & 0 & X_{31}-X_{32}-X_{34}+X_{35} \end{pmatrix} \;,\; J_{mn}^2 = \frac{1}{2} \begin{pmatrix} 0 & 0 & X_{11}-X_{14}-X_{13}+X_{16} \\ 0 & 0 & X_{21}-X_{24}-X_{23}+X_{26} \\ 0 & 0 & X_{31}-X_{34}-X_{33}+X_{36} \end{pmatrix} \quad (14)$$

$$J_{mn}^3 = \frac{1}{2}\begin{pmatrix} 0 & 0 & X_{11} - X_{14} - X_{13} + X_{16} \\ 0 & 0 & X_{21} - X_{24} - X_{23} + X_{26} \\ 0 & 0 & X_{31} - X_{34} - X_{33} + X_{36} \end{pmatrix} \quad (15)$$

Following the above, determinant formula (5), and symbolic manipulations, jacobian matrix determinant J is represented by

$$\begin{aligned} J &= J_0 + \xi J_1 + \eta J_2 + \zeta J_3 + \xi\zeta J_4 + \eta\zeta J_5 + \zeta\zeta J_6 \\ J_0(X_{mi}) &= J(0,0,0) \quad , \quad (m=1,2,3, i=1,..,6) \\ J_1(X_{mi}) &= \frac{\partial J(0,0,0)}{\partial \xi} \; , \; J_2(X_{mi}) = \frac{\partial J(0,0,0)}{\partial \eta} \; , \; J_3(X_{mi}) = \frac{\partial J(0,0,0)}{\partial \zeta} \\ J_4(X_{mi}) &= \frac{\partial J(0,0,0)}{\partial \xi} \; , \; J_5(X_{mi}) = \frac{\partial J(0,0,0)}{\partial \xi} \; , \; J_6(X_{mi}) = \frac{\partial J(0,0,0)}{\partial \xi} \end{aligned} \quad (16)$$

Metric is a function of nodal components and local coordinates. Partial derivatives $J_k$ $(k = 0,..,6)$ are rather lengthy terms. Herein, the use of these precomputed terms is avoided, instead the next interpolations are suggested. The poorest (zero order) approximation for the metric is constant metric - CM.

$$J \approx J_{CM} = J_{cent} \; , \; J_{cent} = \left|J_{mn}^{cent}\right| = \left|J_{mn}(\xi = \frac{1}{3}, \eta = \frac{1}{3}, \zeta = 0)\right| \quad (17)$$

The next (non-unique) interpolation approximates first order or linear metric - LM

$$J \approx J_{LM} = (\frac{11}{12} - \xi - \eta - \zeta)J_1^{LM} + (-\frac{1}{12} + \xi)J_2^{LM} + (-\frac{1}{12} + \eta)J_3^{LM} + (\frac{1}{4} + \zeta)J_4^{LM} \quad (18)$$

$$\begin{aligned} p_1 &: (\xi = \frac{1}{12}, \eta = \frac{1}{12}, \zeta = -\frac{1}{4}) \quad , \quad p_2 : (\xi = \frac{13}{12}, \eta = \frac{1}{12}, \zeta = -\frac{1}{4}) \\ p_3 &: (\xi = \frac{1}{12}, \eta = \frac{13}{12}, \zeta = -\frac{1}{4}) \quad , \quad p_4 : (\xi = \frac{1}{12}, \eta = \frac{1}{12}, \zeta = \frac{3}{4}) \\ J_k^{LM} &= \left|J_{mn}(p_k)\right| \quad , \quad (m,n=1,2,3, k=1,2,3,4) \end{aligned} \quad (19)$$

The above uses metric evaluation at four points $p_k$ $(k = 1,..,4)$. Linear interpolation (18) is not unique. Herein we propose this specific interpolation even though optimality is not rigorously proven. Isoparametric metric - IM approximation $J^{IM} = \varphi^i J_i^{IM}$ $(i = 1,..,6)$, $J_i^{IM} = \left|J(\text{node}\#i)\right|$ yield poorer than (18) accuracy thus omitted from present report.

We develop an *exact* interpolation based on seven points

$$\begin{aligned} J = J_{EX} = \frac{1}{2}((&-\xi - \eta - \zeta + \xi\zeta + \eta\zeta + \zeta^2)J_1^{EX} + (\xi - \xi\zeta)J_2^{EX} + (\eta - \eta\zeta)J_3^{EX} - \\ (&\xi + \eta - \zeta + \xi\zeta + \eta\zeta - \zeta^2)J_4^{EX} + (\xi + \xi\zeta)J_5^{EX} + (\eta + \eta\zeta)J_6^{EX} + 2(1-\zeta^2)J_7^{EX}) \end{aligned} \quad (20)$$

$$J_k^{EX} = |J_{mn}(\text{node}\#k)|, \ (k=1,..,6), \ J_7^{EX} = \left|J_{mn}(\xi=\frac{1}{3},\eta=\frac{1}{3},\zeta=0)\right| \tag{21}$$

The above is identically equal to exact metric (16).

## 4. Results.

CM approximation (17) together with analytical integration (7) result in

$$M_{CM}^{ij} = \rho_0 \frac{J_{cent}}{72} M_{cent}^{ij}, \ (i,j=1,..,6)$$

$$M_{cent}^{ij} = \begin{pmatrix} 4 & 2 & 2 & 2 & 1 & 1 \\ 2 & 4 & 2 & 1 & 2 & 1 \\ 2 & 2 & 4 & 1 & 1 & 2 \\ 2 & 1 & 1 & 4 & 2 & 2 \\ 1 & 2 & 1 & 2 & 4 & 2 \\ 1 & 1 & 2 & 2 & 2 & 4 \end{pmatrix} \tag{22}$$

For lumped mass matrix (8), CM approximation and analytical integration yield

$$M_{CM\,diag}^{ii} = \rho_0 \frac{J_{cent}}{6} M_{cent\,diag}^{ii}, \ (i=1,..,6), \ M_{cent\,diag}^{ii} = \begin{pmatrix} 1 & 1 & 1 & 1 & 1 & 1 \end{pmatrix}^T \tag{23}$$

Where $(\cdot)^T$ denote transpose. Following numerical integration representation (11), it is noted that (22)(23) are equivalent to one point scheme. LM approximation and (7) produce

$$M_{LM}^{ij} = \frac{\rho_0}{4320}(J_1^{LM} M_1^{ij} + J_2^{LM} M_2^{ij} + J_3^{LM} M_3^{ij} + J_4^{LM} M_4^{ij}), \ (i,j=1,..,6) \tag{24}$$

$$M_1^{ij} = \begin{pmatrix} 244 & 98 & 98 & 62 & 19 & 19 \\ 98 & 148 & 74 & 19 & 14 & 7 \\ 98 & 74 & 148 & 19 & 7 & 14 \\ 62 & 19 & 19 & 4 & -22 & -22 \\ 19 & 14 & 7 & -22 & -92 & -46 \\ 19 & 7 & 14 & -22 & -46 & -92 \end{pmatrix}, \ M_2^{ij} = \begin{pmatrix} 28 & 38 & 14 & 14 & 19 & 7 \\ 38 & 124 & 38 & 19 & 62 & 19 \\ 14 & 38 & 28 & 7 & 19 & 14 \\ 14 & 19 & 7 & 28 & 38 & 14 \\ 19 & 62 & 19 & 38 & 124 & 38 \\ 7 & 19 & 14 & 14 & 38 & 28 \end{pmatrix} \tag{25}$$

$$M_3^{ij} = \begin{pmatrix} 28 & 14 & 38 & 14 & 7 & 19 \\ 14 & 28 & 38 & 7 & 14 & 19 \\ 38 & 38 & 124 & 19 & 19 & 62 \\ 14 & 7 & 19 & 28 & 14 & 38 \\ 7 & 14 & 19 & 14 & 28 & 38 \\ 19 & 19 & 62 & 38 & 38 & 124 \end{pmatrix}, \quad M_4^{ij} = \begin{pmatrix} -60 & -30 & -30 & 30 & 15 & 15 \\ -30 & -60 & -30 & 15 & 30 & 15 \\ -30 & -30 & -60 & 15 & 15 & 30 \\ 30 & 15 & 15 & 180 & 90 & 90 \\ 15 & 30 & 15 & 90 & 180 & 90 \\ 15 & 15 & 30 & 90 & 90 & 180 \end{pmatrix} \quad (26)$$

LM approximation (18) and lumped mass matrix formulation (8) produce

$$M_{LM\,diag}^{ii} = \frac{\rho_0}{72}(J_1^{LM} M_{1\,diag}^{ii} + J_2^{LM} M_{2\,diag}^{ii} + J_3^{LM} M_{3\,diag}^{ii} + J_4^{LM} M_{4\,diag}^{ii}), \quad (i=1,..,6) \quad (27)$$

$$\begin{aligned} M_{1\,diag}^{ii} &= \begin{pmatrix} 9 & 6 & 6 & 1 & -2 & -2 \end{pmatrix}^T, \quad M_{2\,diag}^{ii} = \begin{pmatrix} 2 & 5 & 2 & 2 & 5 & 2 \end{pmatrix}^T \\ M_{3\,diag}^{ii} &= \begin{pmatrix} 2 & 2 & 5 & 2 & 2 & 5 \end{pmatrix}^T, \quad M_{4\,diag}^{ii} = \begin{pmatrix} -1 & -1 & -1 & 7 & 7 & 7 \end{pmatrix}^T \end{aligned} \quad (28)$$

Once more, in terms of (11), LM mass matrices are similar to four point rules. Exact metric representation (20) together with analytical integration result in exact consistent and lumped mass matrices

$$M_{EX}^{ij} = \frac{\rho_0}{720} \sum_{k=1}^{7} J_k^{EX} M_k^{ij}, \quad (i,j=1,..,6) \quad (29)$$

$$M_1^{ij} = \begin{pmatrix} 6 & 0 & 0 & -2 & -2 & -2 \\ 0 & -6 & -3 & -2 & -6 & -3 \\ 0 & -3 & -6 & -2 & -3 & -6 \\ -2 & -2 & -2 & -6 & -4 & -4 \\ -2 & -6 & -3 & -4 & -10 & -5 \\ -2 & -3 & -6 & -4 & -5 & -10 \end{pmatrix}, \quad M_2^{ij} = \begin{pmatrix} 6 & 6 & 3 & 2 & 2 & 1 \\ 6 & 18 & 6 & 2 & 6 & 2 \\ 3 & 6 & 6 & 1 & 2 & 2 \\ 2 & 2 & 1 & 2 & 2 & 1 \\ 2 & 6 & 2 & 2 & 6 & 2 \\ 1 & 2 & 2 & 1 & 2 & 2 \end{pmatrix} \quad (30)$$

$$M_3^{ij} = \begin{pmatrix} 6 & 3 & 6 & 2 & 1 & 2 \\ 3 & 6 & 6 & 1 & 2 & 2 \\ 6 & 6 & 18 & 2 & 2 & 6 \\ 2 & 1 & 2 & 2 & 1 & 2 \\ 1 & 2 & 2 & 1 & 2 & 2 \\ 2 & 2 & 6 & 2 & 2 & 6 \end{pmatrix}, \quad M_4^{ij} = \begin{pmatrix} -6 & -4 & -4 & -2 & -2 & -2 \\ -4 & -10 & -5 & -2 & -6 & -3 \\ -4 & -5 & -10 & -2 & -3 & -6 \\ -2 & -2 & -2 & 6 & 0 & 0 \\ -2 & -6 & -3 & 0 & -6 & -3 \\ -2 & -3 & -6 & 0 & -3 & -6 \end{pmatrix} \quad (31)$$

$$M_5^{ij} = \begin{pmatrix} 2 & 2 & 1 & 2 & 2 & 1 \\ 2 & 6 & 2 & 2 & 6 & 2 \\ 1 & 2 & 2 & 1 & 2 & 2 \\ 2 & 2 & 1 & 6 & 6 & 3 \\ 2 & 6 & 2 & 6 & 18 & 6 \\ 1 & 2 & 2 & 3 & 6 & 6 \end{pmatrix}, \quad M_6^{ij} = \begin{pmatrix} 2 & 1 & 2 & 2 & 1 & 2 \\ 1 & 2 & 2 & 1 & 2 & 2 \\ 2 & 2 & 6 & 2 & 2 & 6 \\ 2 & 1 & 2 & 6 & 3 & 6 \\ 1 & 2 & 2 & 3 & 6 & 6 \\ 2 & 2 & 6 & 6 & 6 & 18 \end{pmatrix} \quad (32)$$

$$M_7^{ij} = \begin{pmatrix} 24 & 12 & 12 & 16 & 8 & 8 \\ 12 & 24 & 12 & 8 & 16 & 8 \\ 12 & 12 & 24 & 8 & 8 & 16 \\ 16 & 8 & 8 & 24 & 12 & 12 \\ 8 & 16 & 8 & 12 & 24 & 12 \\ 8 & 8 & 16 & 12 & 12 & 24 \end{pmatrix} \quad (33)$$

The exact lumped mass matrix is given by

$$M_{EX\,diag}^{ii} = \frac{\rho_0}{72} \sum_{k=1}^{7} J_k^{EX} M_{k\,diag}^{ii} \quad , (i=1,..,6) \quad (34)$$

$$\begin{aligned}
M_{1\,diag}^{ij} &= \begin{pmatrix} 0 & -2 & -2 & -2 & -3 & -3 \end{pmatrix}^T, & M_{2\,diag}^{ij} &= \begin{pmatrix} 2 & 4 & 2 & 1 & 2 & 1 \end{pmatrix}^T \\
M_{3\,diag}^{ij} &= \begin{pmatrix} 2 & 2 & 4 & 1 & 1 & 2 \end{pmatrix}^T, & M_{4\,diag}^{ij} &= \begin{pmatrix} -2 & -3 & -3 & 0 & -2 & -2 \end{pmatrix}^T \\
M_{5\,diag}^{ij} &= \begin{pmatrix} 1 & 2 & 1 & 2 & 4 & 2 \end{pmatrix}^T, & M_{6\,diag}^{ij} &= \begin{pmatrix} 1 & 1 & 2 & 2 & 2 & 4 \end{pmatrix}^T \\
M_{7\,diag}^{ij} &= \begin{pmatrix} 8 & 8 & 8 & 8 & 8 & 8 \end{pmatrix}^T
\end{aligned} \quad (35)$$

In terms of (11), EX matrices (29)(34) are analogous to seven-points numerical integration scheme.

## 5. Preliminary numerical study.

Preliminary accuracy examination of new schemes is carried out. Currently used two and nine point schemes are compared to our developments in terms of accuracy and amount of required calculations.

We consider the next coarse mesh element family

$$\begin{aligned}
&X_{11}=0, \quad X_{21}=0, \quad X_{31}=-1, \quad X_{12}=1, \quad X_{22}=0, \quad X_{32}=-1 \\
&X_{13}=0, \quad X_{23}=1, \quad X_{33}=-1, \quad X_{14}=0, \quad X_{24}=0, \quad X_{34}=1+\delta \\
&X_{15}=1+\delta, \quad X_{25}=0, \quad X_{35}=1, \quad X_{16}=0, \quad X_{26}=1+\delta, \quad X_{36}=1
\end{aligned} \quad (36)$$

Where $\delta$ represents mesh coarseness (see Fig. 2). Using the above, metric (16) becomes

$$J^* = \frac{1}{8}(8+12\delta+6\delta^2+\delta^3) - \frac{1}{8}(4\delta+2\delta^2)\xi - \frac{1}{8}(4\delta+2\delta^2)\eta + \frac{1}{8}(8\delta+8\delta^2+2\delta^3)\zeta - \frac{1}{8}(2\delta^2\xi\zeta+2\delta^2\eta\zeta) + \frac{1}{8}(2\delta^2+\delta^3)\zeta^2 \quad (37)$$

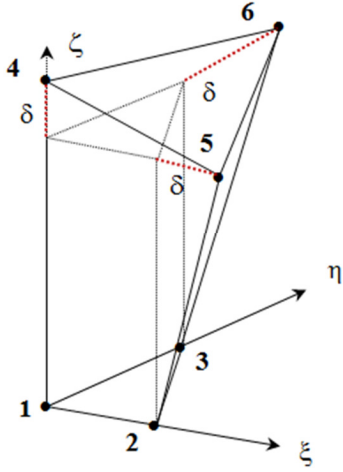

Figure 2: presents coarse element (36), $\delta$ is a deviation from constant metric (parent) element.

Figure 2 and (37) exhibit constant metric element for $\delta = 0$, otherwise coarse mesh element. For each value of coarseness measure $\delta$ (see Fig. 2), consistent CM, LM and EX mass matrix components has been computed, absolute error for each component with respect to exact solution has been evaluated, then averaged among all the components and reported in Figure 3.

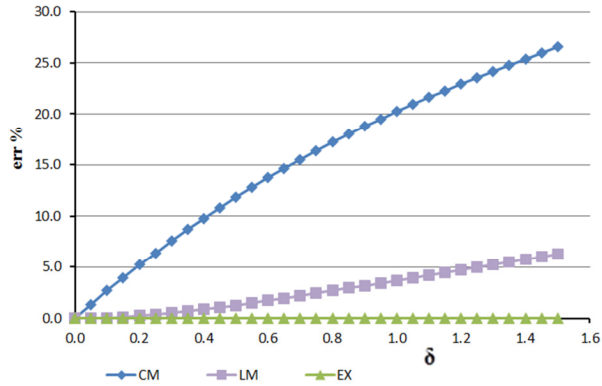

Figure 3: Averaged absolute error of consistent CM, LM and EX mass matrices are reported as a function of the mesh coarseness $\delta$.

CM, LM and EX admit constant metric element $\delta = 0$. LM is considerably more accurate than CM, nevertheless, LM is roughly four times computationally expensive than CM. CM is equivalent to one point scheme (11)(22) while LM is equivalent to four point scheme (11)(24).

Figure 4 records an averaged absolute error of consistent CM mass matrix (22) together with two point numerical integration (11) as a function of mesh coarseness $\delta$.

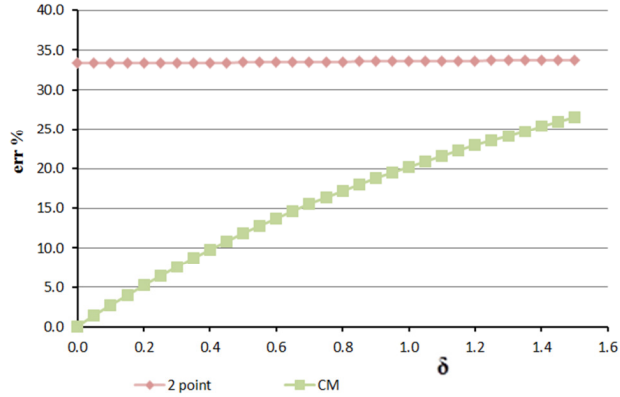

Figure 4: Averaged absolute error of consistent CM rule and mass matrices based on two point numerical integration schemes.

Two-point integration requires roughly twice the amount of computations necessary for CM, yet CM performance is superior throughout all the coarseness range $\delta$. In conclusion, CM based consistent mass matrix is computationally less expensive nevertheless significantly more accurate.

Figure 5 illustrates an averaged absolute error of consistent mass matrix based on nine-point numerical integration and our EX based consistent mass matrix.

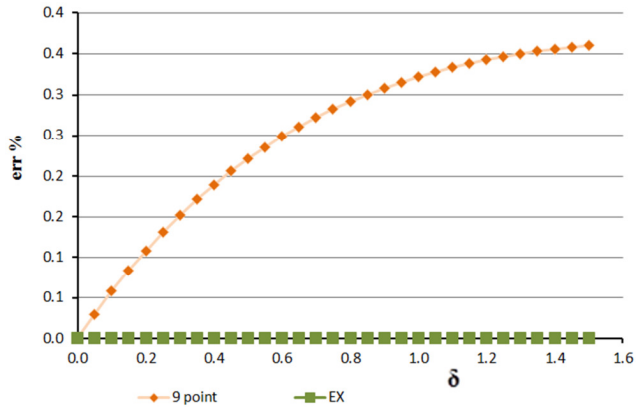

Figure 5: Averaged absolute error of consistent EX based mass matrix vs nine-point numerical integration scheme as a function of the mesh coarseness $\delta$.

EX based mass matrix (29) is equivalent to seven-point numerical scheme thus less expensive than nine point scheme (11), and yet exact. In conclusion, EX based consistent mass matrix is computationally inexpensive but then again exact, therefore superior to currently widely used nine-point rule.

## 6. Conclusions

New integration rules for consistent and lumped mass matrices of six node wedge element (3D solid) have been developed. The metric is approximated as a constant metric (CM), linearly varying metric (LM), special interpolation has been proposed resulting in exact (EX) metric. CM requires one metric evaluation point (the centroid); LM necessitates four evaluation points and linear interpolation functions, while EX uses seven points together with special interpolation functions. Analytical integration resulted in explicit closed-form mass matrices. Coefficient matrices has been defined such that new integration rules take familiar form; CM mass matrix is computationally equivalent to one-point quadrature scheme, LM comparable to four-point numerical integration, and EX mass matrix corresponding to seven-points quadrature rule. Importantly, our EX mass matrix is exact.

Preliminary numerical study has been conducted. Controlled mesh coarseness $\delta$ has allowed gradual coarseness increase while accuracy is examined. Our consistent CM mass matrix exhibits superior accuracy to existing two-point scheme (Fig. 4), although only one evaluation point used. Our consistent LM mass matrix demonstrate considerably more accurate results than CM (Fig. 3), hence LM can be viewed as new 4 point integration scheme for 6 node solid wedge. Contrary to frequently used nine-point quadrature, our EX mass matrix is exact (Fig. 5), however only seven evaluation points are necessary.